\RequirePackage{fix-cm}
\documentclass[smallextended,referee]{svjour3}
\usepackage{latexsym}
\usepackage{amssymb,amsbsy,amsmath,amsfonts,amssymb,amscd}
\usepackage[dvipdfmx]{graphicx}
\usepackage{color}
\usepackage{appendix,comment}

\smartqed
\newcommand{\R}{\mathbb{R}}

\newcommand{\eps}{\varepsilon}
\newcommand{\pt}{\partial}

\title{
\color{black}
Numerical study of the volcano effect in chemotactic aggregation based on a kinetic transport equation with non-instantaneous tumbling
}

\author{Shugo Yasuda}
\institute{Graduate School of Information Science, University of Hyogo, Kobe 650-0047, Japan.\\ \email{yasuda@gsis.u-hyogo.ac.jp}
}
\date{\today}

\begin{document}
\maketitle
\pagestyle{plain}
\pagenumbering{arabic}

\begin{abstract}
Aggregation of chemotactic bacteria under a unimodal distribution of chemical cues was investigated by Monte Carlo (MC) simulation based on a kinetic transport equation, which considers an internal adaptation dynamics as well as a finite tumbling duration.

It was found that there exist two different regimes of the adaptation time, between which the effect of the adaptation time on the aggregation behavior is reversed; that is, when the adaptation time is as small as the running duration, the aggregation becomes increasingly steeper as the adaptation time increases, while, when the adaptation time is as large as the diffusion time of the population density, the aggregation becomes more diffusive as the adaptation time increases. 
Moreover, the aggregation profile becomes bimodal (volcano) at the large adaptation-time regime when the tumbling duration is sufficiently large while it is always unimodal at the small adaptation-time regime.

{\color{black}
A remarkable result of this study is the identification of the parameter regime and scaling for the volcano effect.
That is, by comparing the results of MC simulations to the continuum-limit models obtained at each of the small and large adaptation-time scalings, it is clarified that the volcano effect arises due to the coupling of diffusion, adaptation, and finite tumbling duration, which occurs at the large adaptation-time scaling.
}
\end{abstract}

%
\section{Introduction}
\label{sec:intro}
Chemotactic bacteria, such as \textit{Escherichia coli}, migrate by alternate running and tumbling, where the bacteria run straightly for a certain duration (typically approximately 1 second) and subsequently tumble very quickly (typically, the tumbling duration is approximately 0.1 seconds)~\cite{BB1972}.
Since the length of the run is modulated according to the temporal variation of extracellular chemical cues sensed by the bacteria along their moving pathway, collective dynamics such as traveling waves and aggregations occurs at the population level when the spatial gradients of external chemical cues are sufficiently large.

Kinetic transport equations have been proposed to model the run-and-tumble motion of chemotactic bacteria~\cite{ODA1998,HO2000,EO2004,CMPS2004,DS2005} and successfully utilized to elucidate the mathematics and physics behind the complicated collective dynamics~\cite{SCBPBS2011,EGBAV2016,XXT2018,C2020}.
Since the duration of tumbling is much shorter than the running duration, the tumbling duration is usually ignored, and instead, the instantaneous velocity jump process is considered in the kinetic transport equations.
However, it is fundamental and important to know how that small but finite duration of tumbling affects the collective dynamics.

Recently, studies on non-instantaneous interactions in biological phenomena have become increasingly prevalent.
For example,
{\color{black} in Ref.~\cite{XXT2018}, concentric stripe patterns created by engineering \textit{E.~coli}, which have longer tumbling duration than the wild-type \textit{E.~coli}, was investigated based on a kinetic transport model with non-instantaneous tumbling process.
}
Additionally, in Ref.~\cite{LTLTH2022}, the sub-diffusion of hydration water molecules that stay attached to proteins between jumps was investigated. 
Moreover, in Ref.~\cite{KST2022}, a novel kinetic transport model with non-instantaneous collision operator was proposed to consider the contact inhibition of movement upon cell-cell collisions or collisions under cell-cell adhesion.
The present study also aims to contribute to a fundamental understanding of the effects of the non-instantaneous interaction (i.e., the non-instantaneous tumbling in this study) on a simple chemotactic aggregation of bacteria.

In this paper, we regard the kinetic transport equation, which considers the finite tumbling duration as well as the internal adaptation dynamics~\cite{EO2004}, and clarify their effects on chemotactic aggregation via numerical simulations.
{\color{black}
Especially, we are concerned with the volcano effect (i.e., bimodal aggregation of chemotactic bacteria, which was first observed in an micro-scale experiment~\cite{MBBO2003}).
The volcano effect was previously investigated numerically based on a similar kinetic transport model \cite{SM2011} and an individual-based simulation \cite{JJP2018}.
However, the biological and physical conditions under which the volcano effect takes place and the multiscale mechanisms behind the volcano effect are still largely unknown.
We tackle with these problems by means of a Monte Carlo (MC) simulation of the kinetic transport model.
The results of MC simulations are also compared with the continuum-limit equations obtained at different scalings of the adaptation time, which were previously derived in literatures, e.g., Refs.~\cite{CMPS2004,EO2004,HO2000,PTV2016,PSTY2020,XXT2018}.
An important result of this study is the identification of the parameter regime and scaling for the volcano effect to arise.
That is, it will be clarified that the volcano effect arises due to the coupling of diffusion, adaptation, and finite tumbling duration, which occurs at a certain large adaptation-time scaling.
These results contribute to advance the elucidation of the mathematics behind the volcano effect.

Incidentally, the effect of the adaptation time on the chemotactic aggregation was investigated in our previous study~\cite{Y2021} based on the kinetic transport equation without the finite tumbling duration.
In the previous study, the volcano effect was never observed, but intead, the trapezoidal aggregation, where the aggregation profile in the central region is rather flat, was observed at the large adaptation-time regime.
The present study is an extension of the previous study to include the finite tumbling duration.
Interestingly, as it will be seen in Sec.~\ref{sec:asymp}, although the modification of the continuum-limit equation due to the finite tumbling duration is very small, it is crucial to produce the volcano effect.
}

The rest of the paper is organized as follows:
In Sec.~\ref{sec:formul}, we present the basic kinetic transport model, which considers both the internal adaptation dynamics and the tumbling process. The nondimensionalization of the kinetic transport model is also given.
In Sec.~\ref{sec:asymp}, 
{\color{black}
asymptotic equations of the kinetic transport model, which are utilized to compare with the MC results obtained in Sec.~\ref{sec:numeric}, are summarized.
}
In Sec.~\ref{sec:numeric}, numerical simulations of chemotactic aggregation under an exponential distribution of chemical cues are implemented \textcolor{black}{for one- and two-dimensional spaces} by using the MC code of the kinetic transport model.
From the numerical results, the emergence of bimodal aggregation and its microscopic and macroscopic mechanism are discussed in detail.
Finally, we provide concluding remarks in Sec.~\ref{sec:summary}.

\section{Forumulation}\label{sec:formul}
We consider the run-and-tumble bacteria, where the bacteria modulate their run length according to the memory of the external chemical cue sensed along their moving pathway, while in the tumbling phase, they alter their moving direction in a uniformly random way in a short period of time.

The memory of the external chemical cue is described by the internal state determined via intracellular chemical signal transduction, which has been extensively studied by many authors (for instance, see Refs~\cite{BL1997,KLBTS2005,TSB2008,JOT2010,H2012}).
The mathematical model of intracellular signal transduction is very complicated in general, but the fundamental property necessary for the chemotactic response is described by the excitation and the adaptation process~\cite{SPO1997}.
In the present paper, we consider the following simple adaptation dynamics of the internal state $m \in R$:
$$
\dot m =\frac{M(S)-m}{\tau_a},
$$
where $\tau_a>0$ is the adaptation time and $M(S)$ denotes the local equilibrium of the internal state according to the local concentration of the external chemical cue $S$.

The bacteria modulate their tumbling frequency according to the deviation of the internal state $m$ from the local equilibrium state $M(S)$.
More specifically, the bacteria temporarily decrease (or increase) the tumbling frequency when the local equilibrium state $M(S)$ is higher (or lower) than the current internal state $m$.
The modulated tumbling frequency returns to the basal state when the internal state approaches the local equilibrium state via the adaptation dynamics.
Thus, the adaptation process allows the bacteria to further respond to a subsequent change in the external chemical cue.

Since the modulation of the tumbling frequency is stiff and bounded~\cite{BSB1983}, we write it as follows:
\begin{equation}\label{def:Lambda}
\Lambda_\delta(M(S)-m)=1-F\left(\frac{M(S)-m}{\delta}\right),
\end{equation}
where $\delta$ represents the stiffness of the chemotactic response and the response function $F(X)$ has the following property:
$$
F(0)=0, F'(X)>0, F(x\rightarrow \pm\infty)\rightarrow \pm\chi,
$$
where $\chi$ ($0<\chi<1$) represents the modulation amplitude.

When we write the population density of running cells with a velocity $v \in [v,v+dv]$ and an internal state $m \in [m,m+dm]$ at time $t>0$ and at space $x\in R^d$ \textcolor{black}{(where $d$ is the dimension of space)} as $d\rho_f=f(t,x,v,m)dvdm$ and that of the tumbling cells with an internal state $m\in[m,m+dm]$ as $d\rho_g=g(t,x,m)dm$,
the time evolution of the densities $f(t,x,v,m)$ and $g(t,x,m)$ is described as follows:
\begin{subequations}\label{eq0}
\begin{equation}
    \pt_t f+v\cdot\nabla_x f+
    \pt_m
    \left\{
    \left(\frac{M(S)-m}{\tau_a}
    \right)f
    \right\}
   =
   \mu\frac{g}{||V||} -\lambda\Lambda_\delta(M(S)-m)f
\end{equation}
\begin{equation}
    \pt_t g+\pt_m\left\{\left(\frac{M(S)-m}{\tau_a}\right)g\right\}
    =\lambda\Lambda_\delta(M(S)-m)\int_V f(t,x,v,m)dv-\mu g,
\end{equation}
\end{subequations}
where $\lambda>0$ is the mean tumbling frequency at the basal state, $\mu>0$ is the mean frequency when the tumbling cells change to the running cells, and $||V||$ is the volume of the velocity space, i.e., $||V||=\int_V dv$.
Here, the velocity space is the surface of the ball (i.e., $V=|v|\mathbb{S}^d$).
We also introduce the notation $\nu=\mu^{-1}$, which denotes the mean tumbling duration.

{\color{black}
We also remark that since the duration of tumbling is governed by a Poisson process with a constant rate $\mu$ in Eq.~(\ref{eq0}), very short tumbling events also occur.
Then it is quite an idealization that the post-tumbling direction is uniformly randomized.
The extension to the non-uniform and non-instantaneous tumbling kernel would be an important future study.
}

The total population density $\rho$, the population density of running cells $\rho_f$, and the population density of tumbling cells $\rho_g$ are given as follows:
\begin{subequations}\label{rho0}
\begin{equation}
\rho(t,x)=\rho_f(t,x)+\rho_g(t,x),
\end{equation}
\begin{equation}
    \rho_f(t,x)=\int_{\R}\int_V f(t,x,v,m)dvdm,
\end{equation}
\begin{equation}
    \rho_g(t,x)=\int_{\R}g(t,x,m)dm.
\end{equation}
\end{subequations}

\subsection{Nondimensionalization}
We introduce the nondimensional quantities as follows:
\begin{gather}
    \hat f=f/(\rho_\mathrm{c}/||V||),\quad
    \hat g=g/\rho_\mathrm{c},\quad
    \hat t=t/t_\mathrm{c},\quad
    \hat x=x/L_\mathrm{c},\quad
    \hat v=v/\mathrm{v_c}.
\end{gather}
Here, the subscript ``c'' denotes the characteristic quantities, and we set $\mathrm{v_c}=|v|$ in the following.
Then, we can rewrite (\ref{eq0}) as follows:
\begin{subequations}\label{eq:fg}
\begin{equation}
    \sigma \pt_{\hat t} \hat f+\hat v\cdot\nabla_{\hat x} \hat f+
    \pt_m
    \left\{
    \left(\frac{M(S)-m}{\hat \tau}
    \right)\hat f
    \right\}
   =
   \frac{1}{\eps}\left[\hat \mu \hat g -\Lambda_\delta(M(S)-m)\hat f\right],
\end{equation}
\begin{equation}
    \sigma\pt_{\hat t} \hat g+\pt_m\left\{\left(\frac{M(S)-m}{\hat \tau}\right)\hat g\right\}
    =\frac{1}{\eps}\left[
    \Lambda_\delta(M(S)-m)<\hat f>-\hat \mu \hat g
    \right],
\end{equation}
\end{subequations}
where $<\hat f>$ is the average of $\hat f$ over the velocity space $\hat V$, which is defined as follows:
\begin{equation}
    <f>=\frac{1}{||\hat V||}\int_{\hat V}f(\hat t, \hat x, \hat v, m)d\hat v,
\end{equation}
with $||\hat V||=\int_{\hat V}d\hat v$.
Here, we also introduce the following nondimensional parameters:
\begin{gather}\label{eq:nondim}
    \sigma=L_\mathrm{c}/(\mathrm{v}_\mathrm{c} t_\mathrm{c}),\quad \eps=\mathrm{v}_\mathrm{c}/(\lambda L_\mathrm{c}),\quad
    \hat \tau=\tau_a/(L_\mathrm{c}/\mathrm{v}_\mathrm{c}),\quad \hat \mu=\mu/\lambda.
\end{gather}
Here, $\sigma$ is the time parameter and $\eps$ represents the mean run length at the reference state.

The population densities defined in Eq.~(\ref{rho0}) are written as follows:
\begin{subequations}\label{eq:rho}
\begin{equation}
\hat \rho(\hat t,\hat x)=\hat \rho_f(\hat t,\hat x)+\hat \rho_g(\hat t,\hat x),
\end{equation}
\begin{equation}
    \hat \rho_f(\hat t,\hat x)=\int_{\R}<\hat f>(\hat t,\hat x,m)dm,
\end{equation}
\begin{equation}
    \hat \rho_g(\hat t,\hat x)=\int_{\R}\hat g(\hat t,\hat x,m)dm.
\end{equation}
\end{subequations}

We note that when the tumbling duration is negligibly small compared to the running duration (i.e., $\hat\nu=\hat \mu^{-1}\ll 1$), the density of tumbling cells $\hat g$ becomes negligibly small $\hat g\ll 1$, and thus, Eq.~(\ref{eq:fg}) is reduced as follows:
\begin{equation}\label{eq:kinet_notumble}
\sigma \pt_t \hat f+v\cdot\nabla_x \hat f+\pt_m\left\{
\left(\frac{M(S)-m}{\tau}\right)\hat f
\right\}
=\frac1\eps\Lambda_\delta(M(S)-m)\left(<\hat f>-\hat f\right ).
\end{equation}
This equation was used in a previous study~\cite{Y2021} where instantaneous tumbling events were considered.
Thus, this study is an extension of the previous study to consider the non-instantaneous tumbling events in the chemotactic aggregations.


\section{\textcolor{black}{Continuum-limit equations}}\label{sec:asymp}
\textcolor{black}{It has been proved that} different types of continuum-limit (i.e., $\eps\rightarrow 0$) equations are obtained by the asymptotic analysis of the kinetic transport equation at different scalings of the adaptation time~\cite{EO2004,PTV2016,XXT2018,PSTY2020}.
\textcolor{black}{In this study, we utilize the continuum-limit equations to confirm the asymptotic behaviors of the MC simulations of the kinetic transport equation (\ref{eq:fg}), which will be given in the next section.
In this section, we summarize the continuum-limit equations obtained at two different scalings of the adaptation time.
An asymptotic relation between the two different continuum-limit equations is also briefly explained.}

{\color{black}
We consider small and large adaptation-time scalings, i.e., $\hat \tau=O(\eps)$ and $\hat \tau=O(1/\eps)$, respectively.
We also consider the diffusive time scale $\sigma=\eps$.
}
These settings of the time scale parameters are physically interpreted as follows:
The time scale parameter $\sigma=\eps$ reads that the characteristic time $t_\mathrm{c}$ corresponds to the diffusion time of the population density, i.e., $t_\mathrm{c}=t_d$, where the diffusion time $t_d$ is defined as 
\begin{equation}\label{eq:td}
t_d=L_\mathrm{c}^2/D_\rho
\end{equation}
with the diffusion constant defined as $D_\rho=\mathrm{v_c}^2/\lambda$.
The small adaptation time scaling $\hat \tau=O(\eps)$ indicates that the adaptation time is comparable to the running duration (i.e., $\tau_a\sim \lambda^{-1}$), while the large adaptation-time scaling $\hat \tau=O(1/\eps)$ indicates that the adaptation time is comparable to the diffusion time (i.e., $\tau_a\sim t_d$).

In the following, we only consider the case where the stiffness of the modulation function is moderate, i.e., $\delta=O(1)$ in Eq.~(\ref{def:Lambda}) (although the stiff chemotactic response, such as those considered in the previous studies~\cite{PTV2016,PY2018,PSTY2020} are more realistic in general).
The asymptotic analysis for the stiff chemotactic response with the finite tumbling duration should be an important future work.

Hereafter, we write the nondimensional quantities without "$\hat{\quad}$" for simplicity unless otherwise stated.

\subsection{Small adaptation-time scaling}
We consider the small adaptation-time scaling in Eq.~(\ref{eq:fg}) as follows:
\begin{equation}\label{eq:smalltau}
\tau=\alpha\eps,\quad \sigma=\eps,
\end{equation}
where the parameter $\alpha=O(1)$ denotes the ratio of the adaptation time to the mean running duration (i.e., $\alpha=\tau_a/\lambda^{-1}$).

Then, Eq.~(\ref{eq:fg}) is written as follows:
\begin{subequations}\label{eq:fgalpha}
\begin{equation}
    \eps^2\pt_t f_\eps
    +\eps v\cdot\nabla_xf_\eps
    +\pt_m\left\{
    \left(
    \frac{M_\eps-m}{\alpha}
    \right)f_\eps
    \right\}
    =\mu g_\eps-\Lambda(M_\eps-m)f_\eps,
\end{equation}
\begin{equation}
    \eps^2 \pt_tg_\eps
    +\textcolor{black}{\pt_m}\left\{
    \left(
    \frac{M_\eps-m}{\alpha}
    \right)g_\eps
    \right\}
    =\Lambda(M_\eps-m)<f_\eps>-\mu g_\eps,
\end{equation}
\end{subequations}
where the subscript $\eps$ represents the expansion of the quantity with respect to $\eps$, e.g., $f_\eps=f_0+\eps f_1+\eps^2 f_2 \cdots$.
Here, we also write $M(S_\eps)$ as $M_\eps=M(S_\eps)=M_0+\eps M_1+\cdots$.

The asymptotic analysis of the above equation gives the following standard KS equation for the total population density $\rho$ at the continuum limit $\eps\rightarrow 0$ as follows:
\begin{equation}\label{eq:ModKS}
    \sigma_\nu\pt_t \rho_0-\nabla_x\cdot
    c_d\left[
    \nabla_x \rho_0
    +\frac{\Lambda'(0)\alpha\rho_0}{1+\alpha}\nabla_x M_0
    \right]=0,
\end{equation}
where $\sigma_\nu$ is the time-scale parameter defined as $\sigma_\nu=1+\nu$ and $c_d$ is the diffusion constant calculated as $c_d=1/d$ for $d=1$, 2, and 3.
Here, we note again that the parameter $\nu=\mu^{-1}$ denotes the relative mean tumbling duration to the mean running duration.
The population densities of the running and tumbling cells are obtained as follows:
\begin{gather}
    \rho_f=\frac{1}{1+\nu}\rho_0,
    \quad \rho_g=\frac{\nu}{1+\nu}\rho_0.
\end{gather}
{\color{black}
The formal derivation of the KS equation (\ref{eq:ModKS}) is concisely described in Appendix~\ref{appx:ks}.
}

Equation~(\ref{eq:ModKS}) shows that the tumbling duration $\nu$ only affects the time scale, but the spatial distribution of the population density $\rho_0$ in the steady state are not affected by the tumbling duration at the continuum limit $\eps\rightarrow 0$ in the small adaptation-time scaling (\ref{eq:smalltau}).

\subsection{Large adaptation-time scaling}
We consider the large adaptation time scaling at Eq.~(\ref{eq:fg}) as follows:
\begin{equation}\label{eq:largetau}
\tau=\beta /\eps,\quad \sigma=\eps,
\end{equation}
where the parameter $\beta=O(1)$ denotes the ratio of the adaptation time to the diffusion time (i.e., $\beta=\tau_a/t_d$).

Then, Eq.~(\ref{eq:fg}) is written as follows:
\begin{subequations}\label{eq:fege}
\begin{equation}
    \eps^2\pt_t f_\eps
    +\eps v\cdot\nabla_xf_\eps
    +\eps^2 \pt_m\left\{
    \left(
    \frac{M_\eps-m}{\beta}
    \right)f_\eps
    \right\}
    =\mu g_\eps-\Lambda(M_\eps-m)f_\eps,
\end{equation}
\begin{equation}
    \eps^2 \pt_tg_\eps
    +\eps^2\textcolor{black}{\pt_m}\left\{
    \left(
    \frac{M_\eps-m}{\beta}
    \right)g_\eps
    \right\}
    =\Lambda(M_\eps-m)<f_\eps>-\mu g_\eps.
\end{equation}
\end{subequations}

The asymptotic analysis of Eq.~(\ref{eq:fege}) gives the continuum-limit equation at $\eps\rightarrow 0$ as follows:
\begin{equation}\label{eq:ExKS}
    \pt_t h_0-\nabla_x\cdot\left[
    \frac{c_d}{\Lambda(M_0-m)}\nabla_x\left(
    \frac{h_0}{1+\nu \Lambda(M_0-m)}
    \right)
    \right]
    +\pt_m\left[
    \left(\frac{M_0-m}{\beta}\right)h_0
    \right]=0,
\end{equation}
where $h_0$ denotes the density of the cells with internal state $m$ and is defined as follows:
$$
h_0(t,x,m)=<f_0>(t,x,m)+g_0(t,x,m).
$$
\textcolor{black}{
Although the above continuum-limit equation (\ref{eq:ExKS}) was previously derived in Ref.~\cite{XXT2018}, we concisely describe the formal derivation of Eq.~(\ref{eq:ExKS}) in Appendix~\ref{appx:exks} for the completeness of the present paper.
We also remark that Eq.~(\ref{eq:ExKS}) with $\nu=0$ was also derived in Ref.~\cite{Y2021} from the kinetic transport equation without finite tumbling duration (\ref{eq:kinet_notumble}).
However, interestingly, the volcano effect was not observed in the previous study while, as it will be seen in Sec.~\ref{sec:numeric}, the volcano effect arises both in MC simulations and numerical results of Eq.~(\ref{eq:ExKS}) with the finite tumbling duration $\nu\ne 0$.
This indicates that the small modification introduced in Eq.~(\ref{eq:ExKS}) with the parameter $\nu$ enables to produce the volcano effect.
}

The total population density of cells $\rho_0(t,x)$ is given by the integration of $h_0$ with respect to the internal state $m$, i.e.,
$$
\rho_0(t,x)=\int_\mathrm{R} h_0(t,x,m)dm.
$$
In the following text, we call Eq.~(\ref{eq:ExKS}) the extended Keller-Segel (ExKS) model because the consistency with the standard Keller-Segel model is confirmed at $\beta\rightarrow 0$, as shown in the next subsection.

\subsection{Consistency between the KS and ExKS models}

It is easily seen that by taking the limit as $\alpha\rightarrow \infty$ in Eq.~(\ref{eq:ModKS}), the KS equation is written as follows:
$$
\sigma_\mu\pt_t \rho_0-\nabla_x\cdot c_d\left[
\nabla_x \rho_0+\Lambda'(0)\rho_0\nabla_x M_0
\right]=0.
$$
The same equation is also obtained by taking the limit $\beta\rightarrow 0$ in Eq.~(\ref{eq:ExKS}).
This can be seen as follows:

When taking the limit $\beta\rightarrow 0$ at Eq.~(\ref{eq:ExKS}), we have the equation as follows:
$$
\pt_m\left[(M_0-m)h_0\right]=0.
$$
Thus, the solution $h_0$ at $\beta\rightarrow 0$ is written as follows:
\begin{equation}\label{eq:h0_ks}
h_0(t,x,m)=\rho_0(t,x)\delta(M_0-m),
\end{equation}
where $\delta(m)$ is the Dirac delta function.

On the other hand, by integrating Eq.~(\ref{eq:ExKS}) with respect to $m$, we obtain the equation as follows:
\begin{equation}\label{eq:rhotj1}
    \pt_t \rho_0-\nabla_x\cdot j_1=0,
\end{equation}
where
\begin{equation}\label{def:j1}
    j_1=c_d\int_{\R} B_1(t,x,m)dm,
\end{equation}
and $B_1$ is defined as Eq.~(\ref{eq:B1_h0}).

By substituting Eq.~(\ref{eq:h0_ks}) into Eq.~(\ref{eq:B1_h0}), we obtain the equation as follows:
\begin{align}\label{eq:B1_ks}
    B_1&=-\frac{\mu}{\Lambda(M_0-m)}\nabla_x
    \left(
    \frac{\rho_0\delta(M_0-m)}{\mu + \Lambda(M_0-m)}
    \right)
    \nonumber\\
    &=-\frac{\mu \delta(M_0-m)\nabla_x\rho_0}{\Lambda(M_0-m)(\mu+\Lambda(M_0-m))}
    -\frac{\mu \rho_0\nabla_x M_0}{\Lambda(M_0-m)}
    \left[
    \frac{\delta'(M_0-m)}{\mu+\Lambda(M_0-m)}
    -\frac{\Lambda'(M_0-m)\delta(M_0-m)}{(\mu+\Lambda(M_0-m))^2}
    \right].
\end{align}
Thus, the flux $j_1$ (Eq.~(\ref{def:j1})) at $\beta\rightarrow 0$ is written as follows:
\begin{align}\label{eq:j1}
    j_1&=-\frac{\mu c_d \nabla_x \rho_0}{\Lambda(0)(\mu+\Lambda(0))}
    -\mu c_d\rho_0\nabla_x M_{\textcolor{black}{0}}
    \left[ 
    \frac{\Lambda'(0)(\mu+\Lambda(0))+\Lambda(0)\Lambda'(0)}{\Lambda^2(0)(\mu+\Lambda(0))^2}
    -\frac{\Lambda'(0)}{\Lambda(0)(\mu+\Lambda(0))^2}
    \right]
    \nonumber\\
    &=-\frac{\mu c_d}{\Lambda(0)(\mu+\Lambda(0))}
    \left[
    \nabla_x\rho_0+\frac{\Lambda'(0)}{\Lambda(0)}\rho_0\nabla_xM_{\textcolor{black}{0}}
    \right].
\end{align}
By using $\Lambda(0)=1$ and $\nu=\mu^{-1}$, it is seen that Eq.~(\ref{eq:rhotj1}) with Eq.~(\ref{eq:j1}) provides the KS model at $\alpha\rightarrow \infty$.
Due to this consistency, we can say that the ExKS model (\ref{eq:ExKS}) is an extension of the standard KS model to involve the effects of the internal adaptation dynamics at the large adaptation-time scaling.

\section{Numerical analysis}\label{sec:numeric}

\subsection{Problem and method}
We consider the aggregation of chemotactic bacteria in the one-dimensional interval $-L/2\le x\le L/2$ \textcolor{black}{and two-dimensional square $(x_1,x_2)\in [-L/2,L/2]\times[-L/2,L/2]$} with the periodic boundary condition under the following distribution of the external chemical cues $S(x)$:
\begin{equation}\label{eq:distS}
S(x)=
{\color{black}
\left\{
\begin{array}{cc}
\exp(-|x|),&(\mathrm{for}\quad d=1),\\
\exp(-r),&(\mathrm{for}\quad d=2),
\end{array}
\right.
}
\end{equation}
\textcolor{black}{where $r=\sqrt{x_1^2+x_2^2}$ is the radial distance in the two-dimensional space $d=2$.}
We also consider the logarithmic sensing of the external chemical cue~\cite{KJTW2009}, where the local equilibrium of the internal state is modeled as $M(S)=\log S$.
Thus, for Eq.~(\ref{eq:distS}), $M(S)$ is described as follows:
\begin{equation}\label{eq:MS}
    M(S)=
    {\color{black}
    \left\{
    \begin{array}{cc}
    -|x|,& (d=1),\\
    -r,& (d=2).
\end{array}
\right.
}
\end{equation}
The modulation of tumbling frequency in the chemotactic response is determined by the deviation of the internal state $m$ from the local equilibrium state as described in Eq.~(\ref{def:Lambda}).
In the following numerical simulations, we consider the response function $F(X)$ in Eq.~(\ref{def:Lambda}), which is described as follows:
\begin{equation}\label{eq:F}
    F\left(X\right)=\frac{\chi X}{\sqrt{1+X^2}}.
\end{equation}

The kinetic transport equation (\ref{eq:fg}) is numerically solved \textcolor{black}{for one- and two-dimensional spaces} by the Monte Carlo method described in Appendix~\ref{appx:mc}, while the ExKS model (\ref{eq:ExKS}) is calculated \textcolor{black}{only for one-dimensional space} by \textcolor{black}{a standard} finite-volume (FV) scheme, which is \textcolor{black}{similar to that} described in \textcolor{black}{Appendix~B.1 in Ref.~\cite{Y2021}}.
In the MC simulations, we set the number of mesh intervals as $I=100$, the total number of MC particles as $N=720000$, and the time step size as $\Delta t=2\times 10^{-4}$
{\color{black}
for one-dimensional problem. For two-dimensional problem, we set the number of lattice cells as $50\times 50$ (i.e., $I=50$ in Appendix~\ref{appx:mc}), the total number of MC particles as $N=1.8\times 10^7$, and the time step size as $\Delta t=2\times 10^{-4}$.
}
Here, it should be noted that in the MC method, the characteristic time \textcolor{black}{$t_\mathrm{c}$} is fixed as \textcolor{black}{$t_\mathrm{c}=L_\mathrm{c}/\mathrm{v_c}$} independent of $\eps$, which reads $\sigma=1$ at Eq.~(\ref{eq:fg}), but the time period of the MC simulation is varied according to $\eps$ as $0\le t\le \frac{2L^2}{\eps}$ in order to compare the MC results with the continuum-limit models obtained in the diffusion time scaling, i.e., $\sigma=\eps$.
The MC results are also time-averaged over the time period $\delta t=\frac{0.1L^2}{\eps}$.
In the FV method of the ExKS model, we set the number of mesh intervals as $I=100$ in space and $K=800$ in the internal state, where the internal state variable $m$ is discretized as $m_k=-Y+k \Delta y$ ($k=0,\cdots, K$) with $\Delta y=2Y/K$ and $Y=5$.
The time-step size is set as $\Delta t=1\times 10^{-4}$, and the simulations are performed over the time period $0<t\le 25$.

Numerical simulations are performed for various values of the adaptation time $\tau$, run length $\eps$, stiffness $\delta$, and modulation amplitude $\chi$, while the spatial extent $L=10$ is fixed.
{\color{black}
In this paper, most of the numerical results are given for the spatially one-dimensional problem.
Only a few results are presented for the two-dimensional problem in order to demonstrate that the volcano effect is certainly reproduced by the kinetic transport equation (\ref{eq:fg}) in two-dimensional space.
}

\subsection{Occurrence of bimodal aggregation}\label{sec:occur}
\begin{figure}[htbp]
    \centering
    \includegraphics[width=0.7\textwidth]{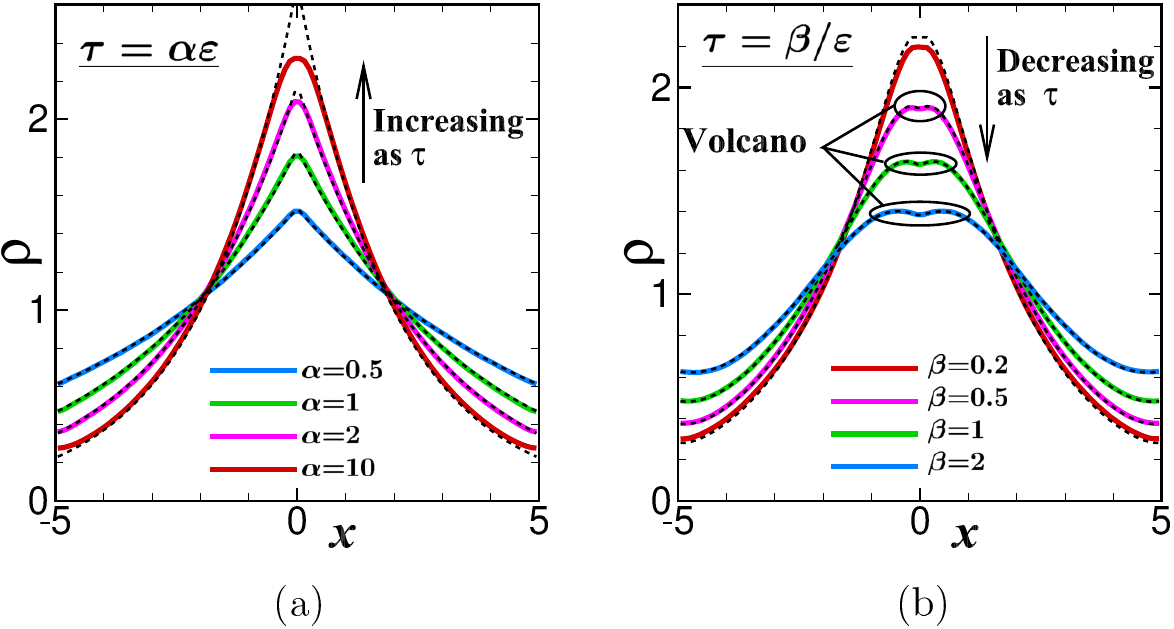}
    \caption{Spatial distributions of the total population density $\rho$ in the steady state at different adaptation times $\tau$ \textcolor{black}{for the one-dimensional problem}.
    Figure (a) shows the results at the small adaptation-time scaling, i.e., $\tau=\alpha \eps$, and Figure (b) shows the results at the large adaptation-time scaling, i.e., $\tau=\beta/\eps$.
    The parameters $\eps=0.1$, $\nu=0.3$, $\delta=1.25$, and $\chi=0.7$ are fixed.
    The colored solid lines show the results of the MC simulations in both figures, while the black dashed lines, which almost overlap with the colored solid lines, show the results of the KS model (\ref{eq:ModKS}) in Fig.~(a) and those of the ExKS model (\ref{eq:ExKS}) in Fig.~(b).
    }
    \label{fig:comptaurho}
\end{figure}

\begin{figure}[htbp]
    \centering
    \includegraphics[width=0.7\textwidth]{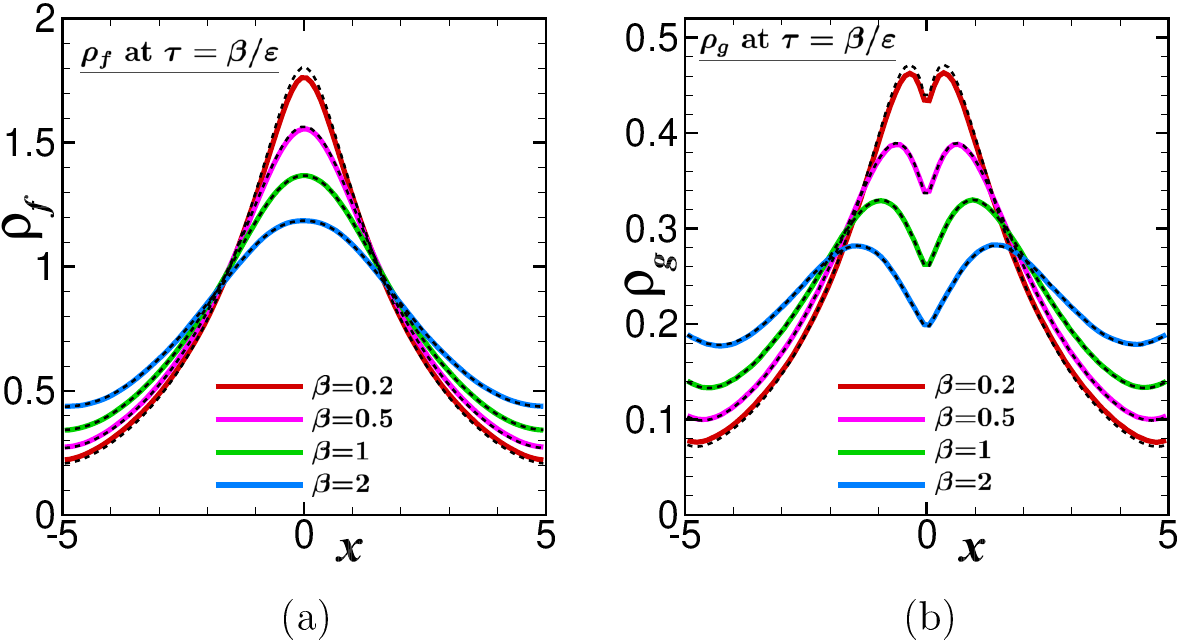}
    \caption{Spatial distributions of the population densities of running and tumbling cells in the steady states, i.e., $\rho_f$ [in (a)] and $\rho_g$ [in (b)] at the large adaptation time scaling $\tau=\beta/\eps$ with $\beta$=0.2, 0.5, 1.0, and 2.0 \textcolor{black}{for the one-dimensional problem}.
    The parameter values of $\eps$, $\nu$, $\delta$, and $\chi$ are the same as those in Fig.~\ref{fig:comptaurho}.
    For the line types, see the caption in Fig.~\ref{fig:comptaurho}.
    }
    \label{fig:comptaurhog}
\end{figure}

Figure \ref{fig:comptaurho} shows the comparison of the aggregation profiles of the total population density $\rho$ at different values of the adaptation time $\tau$.
It is clear that there is a nonmonotonic dependence of the adaptation time on the aggregation profile between the small and the large adaptation-time scalings; at the small adaptation-time scaling $\tau=O(\eps)$ [in Fig~\ref{fig:comptaurho}(a)], the peak of the aggregation profile increases as $\tau$, while at the large adaptation-time scaling $\tau=O(\eps^{-1})$ [in Fig.~\ref{fig:comptaurho}(b)], it decreases as $\tau$.

This nonmonotonic dependence was also observed in a previous study~\cite{Y2021}, where the self-organized aggregation of bacteria without the tumbling phase (i.e., $\nu=0$) was considered.
However, a remarkable difference between the cases with and without the tumbling phase can be observed at the large-adaptation time scaling.
In Fig.~\ref{fig:comptaurho}(b), the volcano effect (i.e., bimodal aggregation) is observed at $\beta=$0.5, 1, and 2.

\begin{figure}[htbp]
    \centering
    \includegraphics[width=0.7\textwidth]{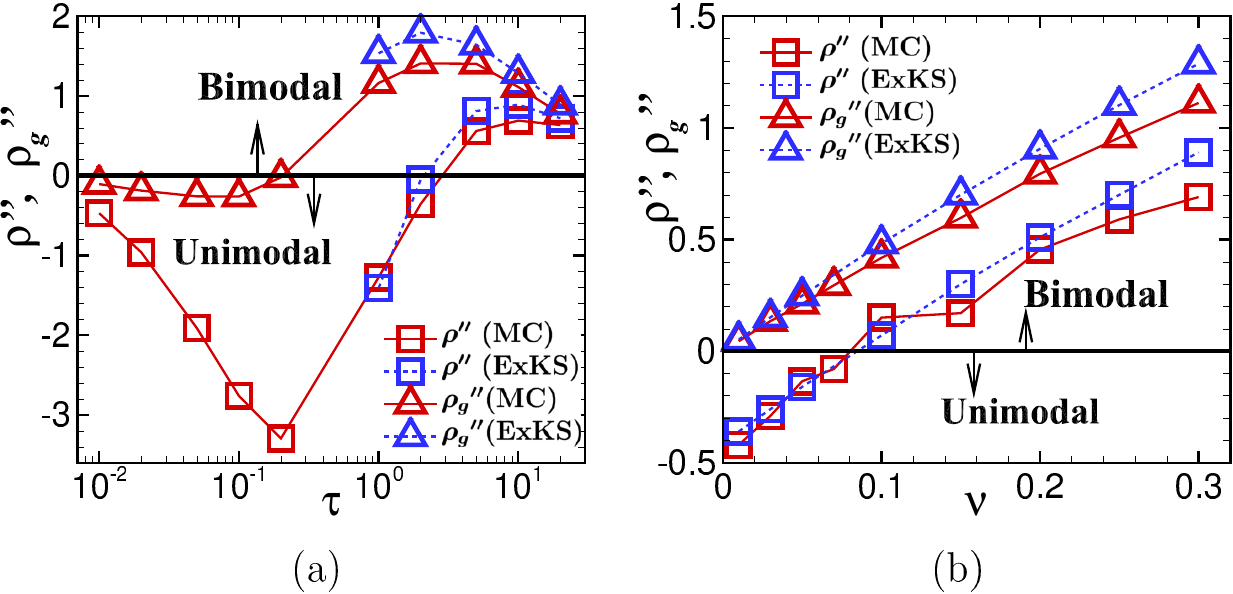}
    \caption{The diagrams of bimodal aggregation with respect to the adaptation time $\tau$ [in (a)] and the tumbling duration $\nu$ [in (b)] \textcolor{black}{for the one-dimensional problem}.
    The second-order derivatives of the total population density $\rho$ and the population density of the tumbling cells $\rho_g$ at $x=0$, $\rho''$ and $\rho''_g$ are shown.
    Thus, $\rho''>0$ and $\rho''_g>0$ represent the occurrence of bimodal aggregation.
    In both figures, the parameters $\eps=0.1$, $\delta=1.25$, and $\chi=0.7$ are commonly fixed, while $\nu=0.3$ is set in (a) and $\tau=10.0$ is set in (b).
    The results obtained by the MC simulation and ExKS model are shown in each figure.
    }
    \label{fig:dderi}
\end{figure}

To clarify the contribution of the tumbling cells to bimodal aggregation, we show the spatial distributions of the running and tumbling cells in Fig.~\ref{fig:comptaurhog}(a) and Fig.~\ref{fig:comptaurhog}(b), respectively.
It is evident that at the large adaptation-time regime, the distribution of running cells $\rho_f$ is unimodal, while that of the tumbling cells $\rho_g$ is bimodal, where the hollow at the central region becomes increasingly larger as $\beta$ increases.
At the small adaptation-time regime, say $\tau\lesssim 0.1$, the distribution of tumbling cells is unimodal, which can be seen from the red-colored triangles in Fig.~\ref{fig:dderi}(a).
Thus, the bimodal distribution of the total population density $\rho$ observed in Fig.~\ref{fig:comptaurho}(b) is due to the local decrease of the tumbling cells at the central region.

Figure~\ref{fig:dderi} is the diagram of the occurrence of bimodal aggregation with respect to the adaptation time $\tau$ and the tumbling duration $\nu$.
Here, the second derivatives of the total population density at $x=0$ are calculated from the numerical results as follows:
$$
\rho''=\frac{\rho_{\frac{I}{2}+1}-\rho_\frac{I}{2}-\rho_{\frac{I}{2}-1}+\rho_{\frac{I}{2}-2}}{(\Delta x)^2},
$$
where $\rho_i$ is the averaged population density in the interval $x\in[x_i,x_{i+1}]$ with $x_i=i\Delta x$ and $x_I=L$.

It can be seen that bimodal aggregation occurs when both the adaptation time $\tau$ and the tumbling duration $\nu$ are sufficiently large, i.e., $\tau\gtrsim 5$ and $\nu\gtrsim 0.1$.
The effect of the tumbling duration $\nu$ on the bimodal aggregation is also seen in \textcolor{black}{SI~3 in the supplemental information (see Online Resource, SI.pdf)}.
It can be seen that the hollow at the central region becomes increasingly deeper as $\nu$ increases.

Notably, it is also seen that the ExKS model (\ref{eq:ExKS}) can reproduce the transient behavior between unimodal and bimodal aggregations described by the kinetic transport equation (although the deviation of the ExKS model from the MC results increases as the tumbling duration increases).
This indicates that the ExKS model (\ref{eq:ExKS}) inherits some essential mathematical structure that should be necessary to describe the volcano effect.

{\color{black}
Figure~\ref{fig:2dmc} shows the MC results for the two-dimensional problem.
It is clearly seen that the volcano effect arises in the two-dimensional space at the large adaptation-time regime.
As in the one-dimensional problem, the population density of running cells is unimodal while that of tumbling cells is bimodal.
The spatial distribution of the local mean run length $\bar \xi$, which is defined as Eq.~(\ref{eq:xi}), sharply increases as the radial distance $r$ approaches to the center of the aggregation.
This feature is also similar to that observed in the local mean run length in the one-dimensional problem, which will be discussed in detail in the next section.
\begin{figure}[t]
\centering
\includegraphics[width=0.9\textwidth]{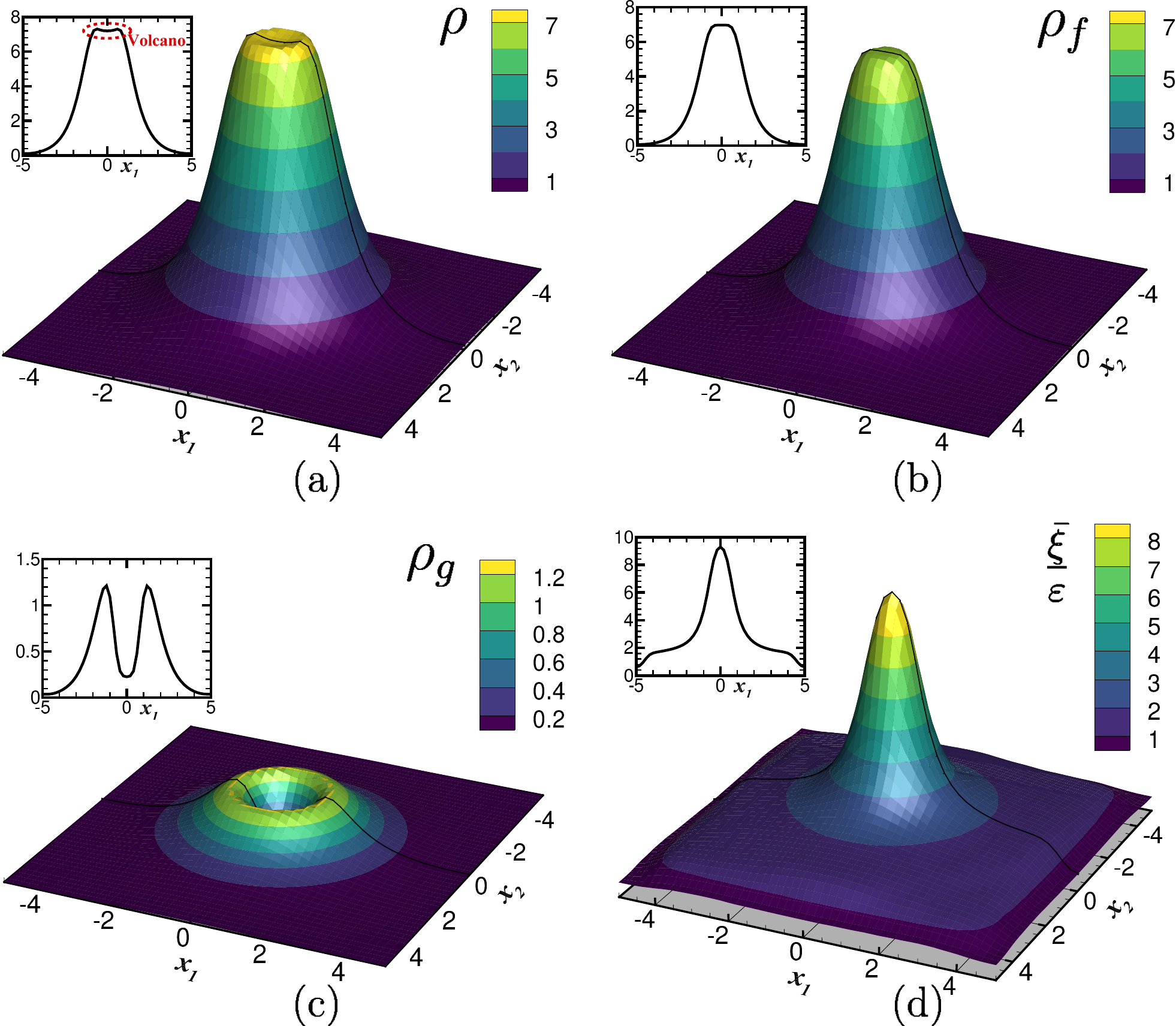}
\caption{
Volcano aggregation in two-dimensional space.
Spatial distributions of the total population density $\rho$ (in (a)), population density of running cells $\rho_f$ (in (b)), population density of tumbling cells $\rho_g$ (in (c)), and local mean run length $\bar \xi$ in the two-dimensional square with $L=10$ are shown.
Each inset show the y-distributions at the intersection $x=5$ (which are shown in sold black lines on the surface of the two-dimensional distributions).
The parameters are set as $\eps=0.1$, $\tau=10$, $\nu=0.3$, $\delta=0.1$, and $\chi=0.9$.
}\label{fig:2dmc}
\end{figure}
}
{\color{black}
We note that MC results for the two-dimensional problem at different values of $\eps$ and $\tau$ are also given in SI~4-6 in the supplemental information (see Online Resource, SI.pdf).
}

\subsection{Why volcano occurs}\label{sec:whyvolcano}
\begin{figure}[htbp]
    \centering
    \includegraphics[width=0.9\textwidth]{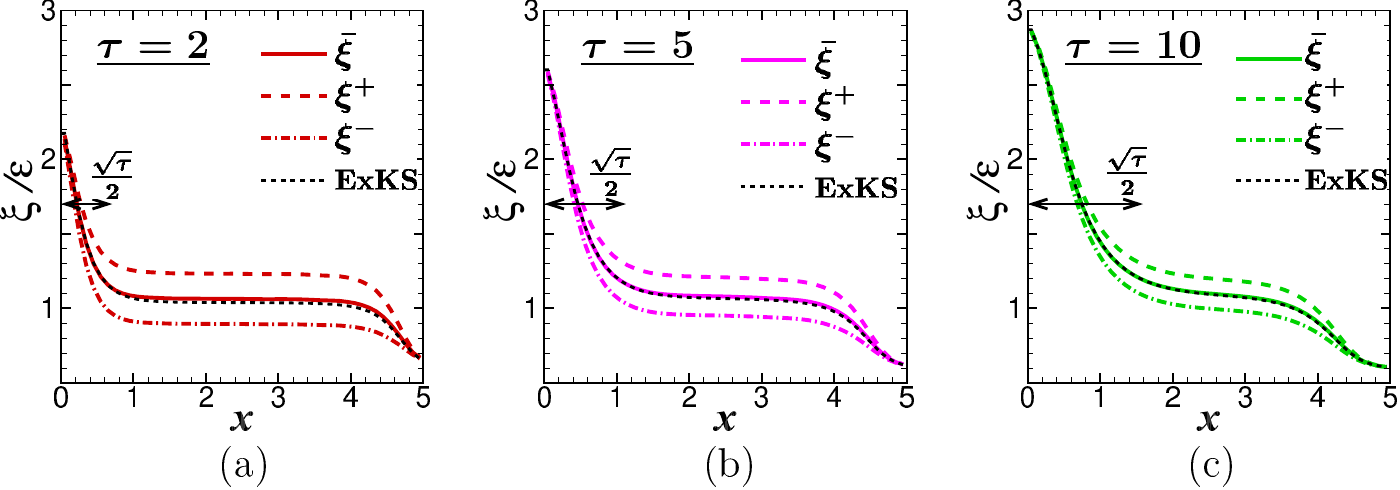}
    \caption{Spatial distributions of the local mean run length defined by Eq.~(\ref{eq:xi}) at $\tau=1$ (a), $\tau=5$ (b) and $\tau=10$ (c).
    The parameters $\eps=0.1$, $\nu=0.3$, $\delta=0.25$, and $\chi=0.7$ are fixed.
    \textcolor{black}{The downward arrows on the horizontal axis in each figure show the position at $x=\sqrt{\eps\tau}$, which represents the thickness of the diffusion layer within the adaptation time.}
    The colored lines show the results obtained by the MC simulations, while the black dashed lines show the results of $\bar \xi$ obtained by the ExKS model (\ref{eq:ExKS}).
    }
    \label{fig:runl}
\end{figure}

In this subsection, we further investigate the microscopic mechanism of the volcano effect \textcolor{black}{with using the numerical results for the one-dimensional problem}.
Figure~\ref{fig:runl} shows the spatial distributions of the local mean run length of the bacteria defined by the equation as follows:
\begin{equation}\label{eq:xi}
    \xi^\pm(t,x)=\int_R \frac{\eps}{\Lambda(M(S)-m)}<f>^\pm(t,x,m) dm,
\end{equation}
where
$$
<f>^\pm=\frac{2}{||V||}\int_{\pm v\cdot\nabla_x S>0} f(v)dv.
$$
Thus, $\xi^\pm$ represent the local mean run lengths of the bacteria, climbing and descending the gradient of chemical cues, respectively, and $\bar \xi=\frac{\xi^++\xi^-}{2}$ is the local mean run length for all of the moving bacteria.

Since the local equilibrium of internal state $M(S)$ is described as Eq.~(\ref{eq:MS}), the spatial gradient of $M(S)$, which is sensed by moving bacteria along their pathway, is constant except the boundaries at $x=0$ and $x=\pm 5$, but it changes stepwise at the boundaries.
The mean run length substantially changes only within the layers near the boundaries with thickness proportional to the diffusion length of the population density during the adaptation time, i.e., $|x|\lesssim \textcolor{black}{\sqrt{\eps\tau}}$ and $|x-5|\lesssim \textcolor{black}{\sqrt{\eps\tau}}$.
We call this layer the diffusion layer since it is created due to the diffusion coupled with the internal adaptation dynamics as explained below.
\textcolor{black}{Here, we note that the non-dimensional form of the diffusion constant defined below Eq.~(\ref{eq:td}) is written as $\hat D_\rho=\eps$.
Thus, $\sqrt{\eps\tau}$ represents the thickness of the diffusion layer within the adaptation time. 
}

In the middle region except the diffusion layers near the boundaries, the bacteria uniformly create a biased random motion according to the local gradient of $M(S)$, where the bacteria climbing the gradient are more likely to have longer run lengths (see $\xi^+$ in Fig.~\ref{fig:runl}), while those descending the gradient are more likely to have shorter run lengths (see $\xi^-$ in Fig.~\ref{fig:runl}).
However, in the close vicinity of the boundary at $x=0$, the bacteria moving across the boundary at $x=0$ with positive velocity are more likely to have longer run lengths than the bacteria moving with positive velocity in the middle region, since they have been climbing the gradient of $M(S)$ in $x<0$.
Thus, the bacteria coming from the left side of the boundary contribute to increasing the local mean run length in the close vicinity of the boundary at $x=0$.

Furthermore, since the bacteria keep the memory of the internal state during the adaptation time $\tau$, which is comparable to the diffusion time in the characteristic length at the large adaptation-time scaling (see the second paragraph in Sec.~\ref{sec:asymp}), the population of the bacteria having a longer run length in the close vicinity at $x=0$ diffuses via the random motions of individual bacteria during the adaptation time.
Thus, the diffusion layer with thickness \textcolor{black}{$\sqrt{\eps\tau}$} is formed near the boundary at $x=0$.
Notably, the increase in the local mean run length in the diffusion layer $|x|<\textcolor{black}{\sqrt{\eps\tau}}$ reduces the local population of the tumbling cells $\rho_g$, as shown in Fig.~\ref{fig:comptaurhog}(b).

The same argument also holds in the right-side diffusion layer; the bacteria moving across the boundary at $x$=5 with a negative velocity are more likely to have shorter run lengths than those of the bacteria moving with negative velocity in the middle region since they have been descending the gradient of the chemical cue in $x>5$ to ensure that they contribute to decreasing the local mean run length.
Thus, as shown in Fig.~\ref{fig:comptaurhog}, the population of tumbling cells slightly increases in the diffusion layer at $|x|=5$.

\begin{figure}[htbp]
    \centering
    \includegraphics[width=0.9\textwidth]{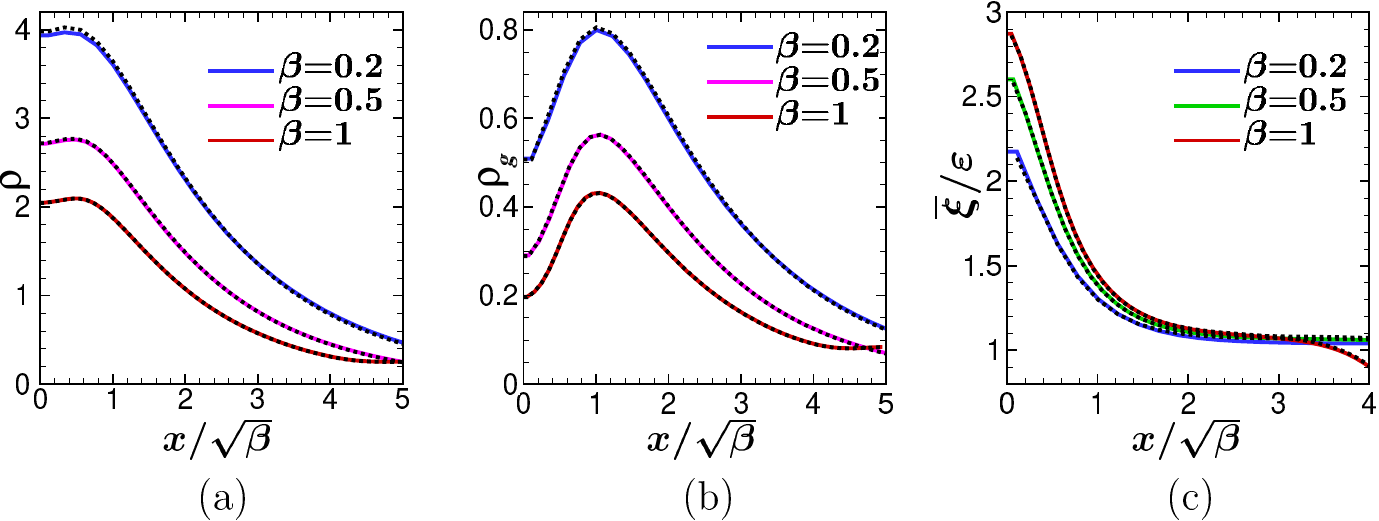}
    \caption{Spatial distributions of the total population density $\rho$ (in (a)), population density of tumbling cells $\rho_g$ (in (b)), and local mean run length $\bar \xi$ in the re-scaled spatial axis $x/\sqrt{\beta}$ at the large adaptation time scaling $\tau=\beta/\eps$.
    The parameters $\nu=0.3$, $\chi=0.7$ and $\delta=0.25$ are fixed.
    The solid lines show the results of the ExKS model, while the dashed line shows the results obtained by MC simulations at $\eps=0.1$.
    }
    \label{fig:ttau}
\end{figure}

The dependence of the adaptation time on the spatial profile is more clearly observed in Fig.~\ref{fig:ttau}, where the spatial distributions of $\rho$, $\rho_g$, and $\bar \xi$ are shown in the scaled coordinate $x/\sqrt{\beta}$.
Notably, the consistency of this scaling property with the ExKS model (\ref{eq:ExKS}) is confirmed; that is, when the equilibrium of the internal state is the linear function described as Eq.~(\ref{eq:MS}), the solution of Eq.~(\ref{eq:ExKS}) has the scaling property as follows:
\begin{equation}
    h_{\delta,\beta}(t,x,m)\propto h_{a\delta,a^2\beta}(a^2 t, a x, am),
\end{equation}
where $a>0$ is an arbitrary constant and $h_{\delta,\beta}(t,x,m)$ is the solution of Eq.~(\ref{eq:ExKS}) with the stiffness $\delta$ and the parameter of the adaptation time $\beta$.
Thus, setting $a=1/\sqrt{\beta}$, we obtain the equation as follows:
\begin{equation}
h_{\delta,\beta}(t,x,m)\propto h_{\frac{\delta}{\sqrt{\beta}},1}\left(\frac{t}{\beta}, \frac{x}{\sqrt{\beta}},\frac{m}{\sqrt{\beta}}\right).
\end{equation}
By integrating the above relation with respect to $m$, we obtain the equation as follows:
\begin{equation}\label{eq:phidel}
\rho_{\delta,\beta}(x)=\rho_{\frac{\delta}{\sqrt{\beta}},1}\left(\frac{x}{\sqrt{\beta}}\right).
\end{equation}
This scaling property indicates that the spatial profile is mainly determined by the parameter $\beta$ when the stiffness $\delta$ does not substantially affect the spatial profile.
Indeed, as seen in \textcolor{black}{SI~1 and SI~2 in the supplemental information (see Online Resource, SI.pdf)}, the variation in the stiffness $\delta$ does not considerably affect the spatial scale of the distribution of population density.
Thus, the scaling property (\ref{eq:phidel}) and the numerical results in SI~1 and SI~2 are consistent with the observation about the dependence of $\beta$ on the spatial profile in Fig.~\ref{fig:ttau}.

\begin{figure}[htbp]
    \centering
    \includegraphics[width=0.9\textwidth]{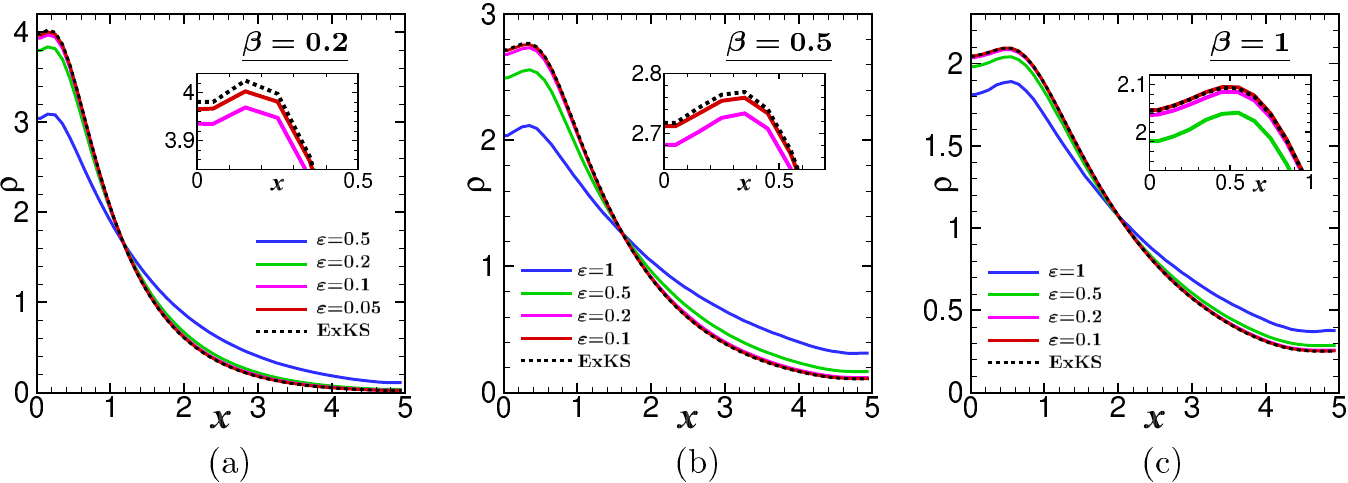}
    \caption{The effect of the run length $\eps$ on the bimodal aggregation and the asymptotic behavior to the continuum limit $\eps\rightarrow 0$.
    The solid lines show the results obtained by the MC simulations, where the adaptation time $\tau$ is set as $\tau=\beta/\eps$ with $\beta=0.2$ in (a), $\beta=0.5$ in (b), and $\beta=1$ in (c), while the black dashed lines show the results of the ExKS model (\ref{eq:ExKS}).
    The parameters $\nu=0.3$, $\chi=0.7$ and $\delta=0.25$ are fixed.
    }
    \label{fig:compeps}
\end{figure}

One may think that the mean run length of the bacteria $\eps$ should directly affect the volcano profile.
Figure~\ref{fig:compeps} shows the effect of the mean run length $\eps$ on the volcano at different values of $\beta$.
The aggregation profile becomes more diffusive as the mean run length $\eps$ increases.
However, the peak position of the aggregation is not affected by the mean run length $\eps$, but it is clearly affected by the parameter $\beta$, as is already shown in Fig.~\ref{fig:ttau} and is explained by the scaling property (\ref{eq:phidel}).
Thus, due to the above numerical results, we conclude that the volcano observed at the large adaptation-time scaling is generated due to the coupling of the internal adaptation and the diffusive motion by the runs and tumbles of bacteria.

It can also be confirmed from the numerical results in Fig.~\ref{fig:compeps} that the ExKS model is a good approximation of the kinetic transport model (\ref{eq:fege}) at the large adaptation-time scaling $\tau=\beta/\eps$ when $\eps$ is moderately small, say $\eps\lesssim 0.2$.


\section{Concluding remarks}\label{sec:summary}
The volcano effect of run-and-tumble chemotactic bacteria was investigated via MC simulations based on the kinetic transport equation (\ref{eq:fg}), which considers a non-instantaneous tumbling as well as an internal adaptation dynamics.
{\color{black}
The MC simulations were performed for one- and two-dimensional problems for a wide range of parameters. 
Especially, in order to identify the parameter regime and scaling for the volcano effect to arise, the results of one-dimensional MC simulations were compared with the numerical results of continuum-limit equations obtained at different scalings of the adaptation time.
}

MC simulations uncovered that the distribution of running cells is always unimodal irrespective to the adaptation time $\tau$, while those of tumbling cells and total population density become bimodal when both the adaptation time $\tau$ and the tumbling duration $\nu$ are sufficiently large, i.e., $\tau=O(1/\eps)$ and $\nu>0.1$.

In order to clarify the microscopic behavior of bacteria forming volcano, the distribution of local mean run length was investigated.
It is clarified that the decrease of the tumbling cells at the central region is caused by the increase of the local mean run length in the central region, where the diffusion layer is created due to the coupling of diffusion and internal adaptation of bacteria.
More concretely, the local mean run length of bacteria increases at the close vicinity of the boundary at $x=0$ due to the mirror symmetry of chemical cues.
The bacteria who have longer run lengths at the central region diffuse due to their individual run-and-tumble motions with keeping the memory of their internal state during the adaptation time.
Thus, the layer in which the local run length is larger than that in the middle region extends to the thickness of the diffusion length in the adaptation time $\sqrt{\eps\tau}$.
The increase of the mean run length indicates the decrease of the tumbling cells, so that the tumbling cells considerably decrease in the diffusion layer and form volcano when the adaptation time is comparable to the diffusion time.
Thus, the coupling of diffusion, adaptation, and tumbling duration is crucial for the formation of volcano.

The fact that the ExKS model (\ref{eq:ExKS}), \textcolor{black}{which is obtained} at the large adaptation-time scaling $\tau=O(1/\eps)$, well approximates the volcano aggregations obtained by MC simulations also convinces that the diffusion of the memory is crucial for the formation of volcano in run-and-tumble bacteria.
The scaling property found in the ExKS model (\ref{eq:phidel}) also clearly shows that the spatial profile of the population density is mainly determined by the diffusion length coupled with the internal adaptation dynamics.

One may think that such a large adaptation time $\tau=O(1/\eps)$ is biologically unrealistic.
However, in the experiment of volcano~\cite{MBBO2003}, the system size is so small as $L_c\sim 100\, \mu\mathrm{m}$, so that the non-dimensional run length $\eps$ is estimated as $\eps=0.2$ by using the typical run length of {\em E. coli} (i.e., $v_c/\lambda$=20 [$\mu$m]).
Furthermore, by using the typical quantities of run duration $\lambda^{-1}$=1 [s], tumbling duration $\mu^{-1}$=0.1 [s], and adaptation time $\tau_a$=10 [s] for \textit{E. coli}, the other non-dimensional parameters defined in Eq.~(\ref{eq:nondim}) are calculated as $\sigma$=0.2, $\tau$=2, and $\nu$=0.1, where $t_c=t_d$ is used in the calculation of $\sigma$.
Thus, the large adaptation-time scaling $\tau=O(1/\eps)$ is not unrealistic but is rather relevant for volcano.

In conclusion we remark that diffusion, adaptation, and non-instantaneous tumbling are key ingredients to describe the non-monotonic aggregation of chemotactic bacteria observed at scales of some tens to hundreds of micro-meters such as volcano, and the kinetic transport model, which considers the non-instantaneous tumbling as well as the internal adaptation dynamics would be useful to elucidate the mathematics behind the complicated aggregation behaviors.


\appendix
\numberwithin{equation}{section}
\setcounter{equation}{0}

\section{Derivation of the continuum-limit model}\label{appx:continuum}
{\color{black}
The continuum-limit models, which are utilized to compared with the MC results in the main text, were previously derived in literatures, e.g., Refs.~\cite{CMPS2004,EO2004,HO2000,PTV2016,PSTY2020,XXT2018}.
In this appendix, we briefly describe the derivation of the models for the completeness of this paper.
}
\subsection{Small adaptation-time scaling}\label{appx:ks}
By following the procedure in Ref.~\cite{PSTY2020}, we change the variables of Eq.~(\ref{eq:fgalpha}) as follows:
$$
f_\eps(t,x,v,m)=p_\eps\left(t,x,v,y=\frac{M(S_\eps)-m}{\eps}\right ),\quad
g_\eps(t,x,m)=q_\eps\left(t,x,y=\frac{M(S_\eps)-m}{\eps}\right).
$$
Then, we have
\begin{subequations}\label{eq:peqe}
\begin{equation}
    \eps^2\pt_t p_\eps +\eps v\cdot\nabla_x p_\eps
    +\pt_y\left\{
    \left(
    v\cdot G_\eps-\frac{y}{\alpha}
    \right)p_\eps
    \right\}
    =\mu q_\eps-\Lambda(\eps y)p_\eps,
\end{equation}
\begin{equation}
    \eps^2\pt_t q_\eps+\pt_y\left\{
    -\frac{y}{\alpha}q_\eps
    \right\}
    =\Lambda(\eps y)<p_\eps>-\mu q_\eps,
\end{equation}
\end{subequations}
where we write $G_\eps=\nabla_x M(S_\eps)$.
By integrating the sum of the above equations with respect to $v$ and $y$, we obtain the following conservation law:
\begin{equation}\label{eq:rho_j}
    \pt_t \rho_\eps +\nabla_x\cdot\left(\frac{j_\eps}{\eps}\right)=0,
\end{equation}
where $\rho_\eps$ is the total population density, i.e., $\rho_\eps=\rho_{p_\eps}+\rho_{q_\eps}$ with $\rho_{p_\eps}=\int <p_\eps>dy$ and $\rho_{q_\eps}=\int q_\eps \textcolor{black}{d}y$, and the flux $j_\eps$ is defined as
\begin{equation}
    j_\eps=\int <v p_\eps>dy.
\end{equation}
As can be seen below, the KS equation is obtained from Eq.~(\ref{eq:rho_j}) at the continuum limit ($\eps\rightarrow 0$).

We assume $p_\eps$ and $q_\eps$ are compactly supported with respect to $y$ and $\Lambda(\eps y)$ in Eq.~(\ref{eq:peqe}) can be expanded as $\Lambda(\eps y)=1+\eps \Lambda'(0) y+O(\eps^2)$.
Then, from the leading-order terms of Eq.~(\ref{eq:peqe}), we obtain the leading-order equation as follows:
\begin{subequations}\label{eq:p0q0}
\begin{equation}\label{eq:p0q0a}
    \pt_y\left\{
    \left(
    v\cdot G_0 -\frac{y}{\alpha}
    \right)p_0
    \right\}=\mu q_0-p_0,
\end{equation}
\begin{equation}
    \pt_y\left\{
    -\frac{y}{\alpha}q_0
    \right\}=<p_0>-\mu q_0.
\end{equation}
\end{subequations}
%
By integrating each of Eqs.~(\ref{eq:p0q0}) w.r.t $y$, we obtain the following relation between the leading-order population densities:
\begin{equation}\label{eq:intp0y}
    \int p_0 dy=\mu \rho_{q_0}=\rho_{p_0}=\frac{\mu}{1+\mu}\rho_0.
\end{equation}
Furthermore, by \textcolor{black}{taking the moment of the above equation multiplied by $v$}, we obtain the flux $j_0$ as follows: 
\begin{equation}
    j_0=\int <vp_0>dy =0.
\end{equation}

From the first-order terms of Eq.~(\ref{eq:peqe}), we have the equation as follows:
\begin{subequations}\label{eq:p1q1}
\begin{equation}
    \nabla_x\cdot (vp_0)
    +\pt_y\left\{
    \left(
    v\cdot G_0-\frac{y}{\alpha}
    \right)p_1
    \right\}
    + \pt_y(v\cdot G_1 p_0)
    =\mu q_1-p_1-\Lambda'(0)yp_0,
\end{equation}
\begin{equation}
    \pt_y\left\{
    \left(
    -\frac{y}{\alpha}
    \right)q_1
    \right\}=<p_1>-\mu q_1+\Lambda'(0)y<p_0>.
\end{equation}
\end{subequations}
%
{\color{black}
By integrating the sum of above equations multiplied by $v$ with respect to $v$ and $y$, we obtain
\begin{equation*}
\begin{split}
\nabla_x\cdot \int<v\otimes v p_0>dy=-j_1-\Lambda'(0)\int y<vp_0>dy,\\
\frac{\mu}{1+\mu}c_d\nabla_x\rho_0=-j_1-\Lambda'(0)\int y<vp_0>dy,
\end{split}
\end{equation*}
where we use Eq.~(\ref{eq:intp0y}) and $<v\otimes v>=c_d I$ with $c_d=1/d$ for $d=1,2,3$.
Here, $I$ is the identity matrix.
}
The last term of the above equation is obtained by integrating Eq.~(\ref{eq:p0q0a}) multiplied by $vy$ with respect to $v$ and $y$:
\begin{equation*}
\begin{split}
-\int y\pt\left\{
<v\otimes v p_0>\cdot G-\frac{y}{\alpha}<vp_0>
\right\}dy
=\int y<vp_0>dy,\\
\int <v\otimes vp_0>dy\cdot G-\frac{1}{\alpha}\int y<vp_0>dy=\int y<vp_0>dy,\\
\int y<vp_0>dy=\frac{\mu}{1+\mu}\frac{\alpha}{1+\alpha}c_d\rho_0 G.
\end{split}
\end{equation*}
Hence, the flux $j_1$ is written as follows:
\begin{equation}
    j_1=-\frac{\mu}{1+\mu}c_d\left[
    \nabla_x \rho_0+\frac{\alpha}{1+\alpha}\Lambda'(0)G\rho_0
    \right].
\end{equation}

Thus, by taking the limit $\eps\rightarrow 0$ at Eq.~(\ref{eq:rho_j}), we obtain the KS equation (\ref{eq:ModKS}).

\subsection{The Extended Keller-Segel model}\label{appx:exks}
The derivation of Eq.~(\ref{eq:ExKS}) is as follows.
In the following, we write the average of $f_\eps$ over the velocity space as $A_\eps=<f_\eps>$.

From the leading-order terms of Eq.~(\ref{eq:fege}), we can write the leading-order solution as follows:
\begin{equation}\label{eq:f0g0}
f_0=A_0(t,x,m),\quad g_0=\frac{\Lambda(M_0-m)}{\mu}A_0(t,x,m).
\end{equation}
Here, $A_0(t,x,m)$ is an unknown function independent of the velocity $v$.

From the $\eps^1$ terms of Eq.~(\ref{eq:fege}), we obtain
\begin{subequations}\label{eq:e1terms}
\begin{equation}
    v\cdot\nabla_x A_0=\mu g_1-\Lambda(M_0-m) f_1-\Lambda'(M_0-m)M_1A_0,
\end{equation}
\begin{equation}
    0=\Lambda(M_0-m)A_1+\Lambda'(M_0-m)M_1A_0-\mu g_1.
\end{equation}
\end{subequations}
By taking the sum of the above equations, we obtain the following equation,
$$
v\cdot\nabla_x A_0=\Lambda(M_0-m)(A_1-f_1).
$$
Hence, $f_1$ can be written in the form
\begin{subequations}\label{eq:f1}
\begin{equation}
f_1=A_1(t,x,m)+v\cdot B_1(t,x,m),    
\end{equation}
with
\begin{equation}\label{eq:B1}
    B_1=-\frac{\nabla_x A_0}{\Lambda(M_0-m)}.
\end{equation}
\end{subequations}
From Eq.~(\ref{eq:e1terms}b), we can also write $g_1$ as follows:
\begin{equation}
    g_1=\frac1\mu\left(\Lambda(M_0-m)A_1+\Lambda'(M_0-m)M_1A_0\right).
\end{equation}

Subsequently, from the $\eps^2$ terms of Eq.~(\ref{eq:fege}), we obtain the equation as follows:
\begin{subequations}\label{eq:e2terms}
\begin{equation}
\begin{split}
    \pt_t A_0+v\cdot \nabla_x (A_1+v\cdot B_1)
    +\pt_m\left\{
    \left(
    \frac{M_0-m}{\beta}
    \right)A_0
    \right\}
    =\mu g_2-\Lambda(M_0-m) f_2\\
    -\Lambda'(M_0-m)M_1(A_1+v\cdot B_1)-\frac12\Lambda''(M_0-m)M_2A_0,
\end{split}
\end{equation}
\begin{equation}
\begin{split}
    \frac{\Lambda(M_0-m)}{\mu}\pt_t A_0
    +\pt_m
    \left\{
    \left(
    \frac{M_0-m}{\beta}
    \right)
    \frac{\Lambda(M_0-m)}{\mu}
    A_0
    \right\}
    =\Lambda(M_0-m)A_2-\mu g_2\\
    +\Lambda'(M_0-m)M_1A_1
    \textcolor{black}{+}\frac12\Lambda''(M_0-m)M_2A_0.
\end{split}
\end{equation}
\end{subequations}
By integrating the sum of the above equations with respect to $v$, we obtain the equation to determine the leading-order solution $A_0$ as follows:
\begin{equation}\label{eq:exA0}
\left(
1+\frac{\Lambda(M_0-m)}{\mu}
\right)\pt_tA_0
+\nabla_x\cdot(c_d B_1)
+\pt_m\left\{
\left(
\frac{M_0-m}{\textcolor{black}{\beta}}
\right)
\left(
1+\frac{\Lambda(M_0-m)}{\mu}
\right)A_0
\right\}=0.
\end{equation}

When we write the total density of cells with an internal state $m$ as $h_0=<f_0>+g_0$, i.e., from Eq.~(\ref{eq:f0g0}),
\begin{equation}\label{def:h0}
    h_0(t,x,m)=\left(1+\frac{\Lambda(M_0-m)}{\mu}\right)A_0(t,x,m),
\end{equation}
we can rewrite Eq.~(\ref{eq:exA0}) as follows:
\begin{subequations}
\begin{equation}\label{eq:h0}
\pt_t h_0+\nabla_x\cdot (c_d B_1)
+\pt_m\left\{\left(\frac{M_0-m}{\textcolor{black}{\beta}}\right)h_0\right\}=0,
\end{equation}
with
\begin{equation}\label{eq:B1_h0}
B_1=-\frac{1}{\Lambda(M_0-m)}\nabla_x\left(\frac{h_0}{1+\frac{\Lambda(M_0-m)}{\mu}}\right).
\end{equation}
\end{subequations}
Thus, we obtain the ExKS model~(\ref{eq:ExKS}).

\section{Monte Carlo method}\label{appx:mc}
We extend the Monte Carlo (MC) method developed in Refs.~\cite{Y2017,PY2018,Y2021} to include the tumbling duration.
In the MC method, we use the time scale $t_0=L_0/V_0$ (which reads $\sigma=1$ in Eq.~(\ref{eq:fg})).
We also use the variable $y=M(S)-m$ instead of the internal state $m$ itself, since this change of variable rewrites the internal adaptation dynamics into the formulation involving the material derivative of the chemical cue along the moving pathway of bacteria, which can be calculated in a straightforward way in the MC code.

By changing the variable as $f(t,x,v,m)=p(t,x,v,y=M(S)-m)$ and $g(t,x,m)=q(t,x,y=M(S)-m)$, Eq.~(\ref{eq:fg}) reads as follows:
\begin{subequations}\label{eq:pq}
\begin{equation}
    \pt_t p+v\cdot\nabla_x p+\pt_y\left\{
    \left(D_t M(S)-\frac{y}{\tau}\right)p
    \right\}=\frac{1}{\eps}\left[\mu q-\Lambda(y)p\right],
\end{equation}
\begin{equation}
    \pt_t q+\pt_y\left\{
    \left(\pt_t M(S)-\frac{y}{\tau}\right)q
    \right\}=\frac{1}{\eps}\left[\Lambda(y)<p>-\mu q\right],
\end{equation}
\end{subequations}
where $D_t$ denotes the material derivative along the moving pathway, i.e., $D_t=\pt_t+v\cdot\nabla_x$.
The $y$-derivative term in Eq.~(\ref{eq:pq}) denotes the internal dynamics of individual cells, where the internal state $y$ changes according to the temporal variation of the external chemical cue sensed by individual cells along their moving pathways, i.e., for the logarithmic sensing $M(S)=\log S$, the equation is as follows:
\begin{equation}\label{eq:doty}
\dot y=\frac{D_t S}{S} -\frac{y}{\tau}.
\end{equation}

\textcolor{black}{For the one-dimensional problem ($d=1$),} the velocity space is composed of three discrete velocities, i.e., $v=\{-1,0,1\}$ and the one-dimensional space $-L/2\le x \le L/2$ is divided into the uniform mesh cells $C_i=[i\Delta x,(i+1)\Delta x]$ ($i=0,\cdots,I-1$), where $\Delta x=L/I$ is the width of the mesh interval.
{\color{black}
For the two-dimensional problem ($d=2$), the velocity space is composed of $v=0$ and $v=(\cos\theta,\sin\theta)$, where the angle $\theta$ is uniformly distributed in $[0,2\pi]$, and the two-dimensional space $[-L/2,L/2]\times [-L/2,L/2]$ is divided into the uniform lattice mesh cells $C_i=[i_1\Delta x,(i_1+1)\Delta x]\times [i_2\Delta x,(i_2+1)\Delta x]$ ($i_1,i_2=0,\cdots,I-1$) with $i=i_1+i_2\times I$.
}
Initially, MC particles are uniformly distributed in each mesh cell with the equilibrium internal state at $y=0$.
The initial velocities of each MC particle are randomly determined from $v=\{-1,+1\}$ 
{\color{black}
for $d=1$ while, for $d=2$, the direction $\theta$ of the individual initial velocities $v=(\cos\theta,\sin\theta)$ are uniformly randomly determined.
}
Thus, the initial condition for Eq.~(\ref{eq:pq}) is described as $p=\delta(y)$ and $q=0$ at $t=0$.

Then, the position $r^k_l$, velocity $v_l^k$, and internal state $y^k_l$ of the $l$th MC particle at time $t=k\Delta t$ ($k\ge 1$) are determined as follows:
\begin{enumerate}
\item Each MC particle moves as follows:
\begin{equation}
r^k_l=r^{k-1}_l+v_l^{k-1}\Delta t.
\end{equation}
\item The population density in the $i$th mesh cell $C_i$, $\rho^k_i$ are calculated as follows:
\begin{equation}
\rho_i^k=\frac{1}{\bar N} \sum_{l=0}^N\int_{C_i}\delta(x-r_l^k)dx,
\end{equation}
where $\bar N$ is the average number of the MC particles in the mesh cell, i.e., $\bar N=N/I$ for $d=1$ and $\bar N=N/I^2$ for $d=2$.

\item The internal state of the $l$th MC particle, $y^k_l$, is updated according to Eq.~(\ref{eq:doty}) as follows:
\begin{equation}\label{eq_yupdate}
\frac{y^k_l-y^{k-1}_l}{\Delta t}=
\frac{S^k_{l}-S^{k-1}_{l}}{\Delta t S^{k-1}_{l}}-\frac{y^k_l}{\tau},
\end{equation}
where $S^k_{l}$ denotes the local concentration of the chemical cue at the position of the $l$th MC particle, i.e., $S^k_{l}=-\exp|r_l^k|$.
We note that in Eq.~(\ref{eq_yupdate}), the pathway derivative $D_t S$ in Eq.~(\ref{eq:doty}) is given by the rate of change of $S$ sensed by each MC particle, i.e., $(S_{l}^k-S_{l}^{k-1})/\Delta t$.
\item The running or tumbling state of the $l$th MC particle is stochastically determined according to the right-hand side of Eq.~(\ref{eq:pq}).
When $v^{k-1}_l\ne 0$, the velocity is changed to $v^k_l=0$ by the probability $\frac{\Delta t \Lambda(y_l^k)}{\eps}$, while when $v^{k-1}_l=0$, the velocity is changed to \textcolor{black}{$|v^{k}_l|=1$} by the probability $\frac{\Delta t \mu}{\eps}$, where the moving direction, \textcolor{black}{i.e., $v=\{-1,+1\}$ for $d=1$ and the angle $\theta\in[0,2\pi]$ of $v=(\cos \theta,\sin\theta)$ for $d=2$}, is uniformly randomly determined.
The particles that are not selected to change their velocities retain their states, i.e., $v_l^k=v_l^{k-1}$.
\item Return to 1.
\end{enumerate}

\begin{acknowledgements}
The author would like to acknowledge partial funding from the Japan-France Integrated action Program PHC SAKURA, Grant No. JPJSBP120193219.
\end{acknowledgements}

%
%


\end{document}